\documentclass[12pt]{article}

\usepackage{amssymb}
\usepackage{amsmath}
\usepackage{graphics}
\usepackage{graphicx}

\newtheorem{theorem}{Theorem}[section]
\newtheorem{definition}[theorem]{Definition}
\newtheorem{proposition}[theorem]{Proposition}

\newtheorem{lemma}[theorem]{Lemma}

\newtheorem{corol}[theorem]{Corollary}

\newenvironment{bew}{\medskip \noindent {\it Proof}:}{\hfill $\Box$\medskip}

\def\R{{\mathbb{R}}}
\def\N{{\mathbb{N}}}
\def\E{{\mathbb{E}}}
\def\P{{\mathbb{P}}}

\def\1{\mbox{{\rm 1 \hskip-1.4ex I}}}

\newcommand{\dint}{{\operatorname{d}}}
\newcommand{\V}{\operatorname{Var}}
\newcommand{\Lip}{\operatorname{Lip}}
\newcommand{\dom}{\operatorname{dom}}

\begin{document}

\title{\textbf{Normal approximation of Poisson functionals in Kolmogorov distance}}
\author{Matthias Schulte\footnote{\textit{Karlsruhe Institute of Technology, Institute of Stochastics, D-76128 Karlsruhe, e-mail: matthias.schulte[at]kit.edu}}}
\date{}
\maketitle

\begin{abstract}
Peccati, Sol\`e, Taqqu, and Utzet recently combined Stein's method and Malliavin calculus to obtain a bound for the Wasserstein distance of a Poisson functional and a Gaussian random variable. Convergence in the Wasserstein distance always implies convergence in the Kolmogorov distance at a possibly weaker rate. But there are many examples of central limit theorems having the same rate for both distances. The aim of this paper is to show this behaviour for a large class of Poisson functionals, namely so-called U-statistics of Poisson point processes. The technique used by Peccati et al.\ is modified to establish a similar bound for the Kolmogorov distance of a Poisson functional and a Gaussian random variable. This bound is evaluated for a U-statistic, and it is shown that the resulting expression is up to a constant the same as it is for the Wasserstein distance.
\end{abstract}
\begin{flushleft}
\end{flushleft}

\section{Introduction and results}

Let $\eta$ be a Poisson point process over a Borel space $(X,{\mathcal{X}})$ with a $\sigma$-finite intensity measure $\mu$ and let $F=F(\eta)$ be a random variable depending on the Poisson point process $\eta$. In the following, we call such random variables Poisson functionals. Moreover, we assume that $F$ is square integrable (we write $F\in L^2(\P_\eta)$) and satisfies $\E F=0$. By $N$ we denote a standard Gaussian random variable. In \cite{Peccatietal2010}, Peccati, Sol\'e, Taqqu, and Utzet derived by a combination of Stein's method and Malliavin calculus the upper bound
\begin{equation}\label{eqn:Peccatietal}
d_W(F,N)\leq \E|1-\langle DF,-DL^{-1}F\rangle_{L^2(\mu)}|+\int_{X}\E(D_zF)^2 \, |D_zL^{-1}F| \, \mu(\dint z)
\end{equation}
for the Wasserstein distance of $F$ and $N$. Here, $\langle\cdot,\cdot\rangle_{L^2(\mu)}$ stands for the inner product in $L^2(\mu)$, and the difference operator $D$ and the inverse of the Ornstein-Uhlenbeck generator $L$ are operators from Malliavin calculus. The underlying idea of these operators is that each square integrable Poisson functional has a representation
$$F=\E F+\sum_{n=1}^\infty I_n(f_n),$$
where the $f_n$ are square integrable functions supported on $X^n$, $I_n$ stands for the $n$-th multiple Wiener-It\^o integral, and the right-hand side converges in $L^2(\P_\eta)$. This decomposition is called Wiener-It\^o chaos expansion, and the Malliavin operators of $F$ are defined via their chaos expansions. The operators $D_zF$ and $D_zL^{-1}F$ that occur in \eqref{eqn:Peccatietal} are given by
$$D_zF=\sum_{n=1}^\infty n \, I_{n-1}(f_n(z,\cdot))\ \text{ and }\ D_zL^{-1}F=-\sum_{n=1}^\infty I_{n-1}(f_n(z,\cdot))\ \text{ for }z\in X.$$
Here, $f_n(z,\cdot)$ stands for the function on $X^{n-1}$ we obtain by taking $z$ as first argument. For exact definitions including the domains and more details on the Malliavin operators we refer to Section \ref{sec:Preliminaries}.

The Wasserstein distance between two random variables $Y$ and $Z$ is defined by
$$d_W(Y,Z)=\sup\limits_{h\in \Lip(1)}|\E h(Y)-\E h(Z)|,$$
where $\Lip(1)$ is the set of all functions $h:\R\rightarrow\R$ with Lipschitz constant less than or equal to one. Another commonly used distance for random variables is the Kolmogorov distance
$$d_K(Y,Z)=\sup\limits_{s\in\R}|\P(Y\leq s)-\P(Z\leq s)|,$$
which is the supremum norm of the difference of the distribution functions of $Y$ and $Z$. Because of this straightforward interpretation, one is often more interested in the Kolmogorov distance than in the Wasserstein distance.
For the important case that $Z$ is a standard Gaussian random variable $N$ it is known (see \cite[Theorem 3.1]{ChenShao}) that
\begin{equation}\label{eqn:dKdW}
d_K(Y,N)\leq 2\sqrt{d_W(Y,N)}.
\end{equation}
This inequality gives us for the Kolmogorov distance a weaker rate of convergence than for the Wasserstein distance. But for many classical central limit theorems, one has actually the same rate of convergence for both metrics.

In order to overcome the problem that a detour around the Wasserstein distance and the inequality \eqref{eqn:dKdW} often gives a suboptimal rate of convergence for the Kolmogorov distance, we derive a similar bound as \eqref{eqn:Peccatietal} for the Kolmogorov distance by a modification of the proof in \cite{Peccatietal2010}.
\begin{theorem}\label{thm:Theorem1}
Let $F\in L^2(\P_\eta)$ with $\E F=0$ be in the domain of $D$ and let $N$ be a standard Gaussian random variable. Then
\begin{eqnarray}\label{eqn:BoundNormalApproximation}
d_K(F,N) &\leq& \E|1-\langle DF,-DL^{-1}F\rangle_{L^2(\mu)}|+2\E\langle (DF)^2,|DL^{-1}F| \rangle_{L^2(\mu)}\\ && + 2\E\langle (DF)^2,|F \ DL^{-1}F| \rangle_{L^2(\mu)} + 2\E\langle (DF)^2,|DF \ DL^{-1}F| \rangle_{L^2(\mu)} \notag\\ && +\sup\limits_{s\in\R}\E\langle D\1(F>s),DF \, |DL^{-1}F|\rangle_{L^2(\mu)}\notag\\
& \leq & \E|1-\langle DF,-DL^{-1}F\rangle_{L^2(\mu)}|+2c(F) \sqrt{\E\langle (DF)^2,(DL^{-1}F)^2\rangle_{L^2(\mu)}} \notag\\
&&+\sup\limits_{s\in\R}\E\langle D\1(F>s),DF \, |DL^{-1}F|\rangle_{L^2(\mu)}\notag
\end{eqnarray}
with
$$c(F)=\sqrt{\E\langle (DF)^2,(DF)^2\rangle_{L^2(\mu)}}+\left(\E\langle DF,DF\rangle_{L^2(\mu)}^{2}\right)^{\frac 1 4}\left(\left(\E F^4\right)^{\frac{1}{4}}+1\right).$$
\end{theorem}
Comparing \eqref{eqn:Peccatietal} and \eqref{eqn:BoundNormalApproximation}, one notes that both terms of the Wasserstein bound \eqref{eqn:Peccatietal} also occur in \eqref{eqn:BoundNormalApproximation}, which means that the bound for the Kolmogorov distance is always larger.

We apply our Theorem \ref{thm:Theorem1} to two situations, where we obtain the same rate of convergence for the Kolmogorov distance and the Wasserstein distance. At first we derive the classical Berry-Esseen inequality with the optimal rate of convergence for the normal approximation of a classical Poisson random variable.
As another application of Theorem \ref{thm:Theorem1}, we consider so-called U-statistics of Poisson point processes, which are defined as
$$F=\sum_{(x_1,\hdots,x_k)\in\eta_{\neq}^k}f(x_1,\hdots,x_k)$$
where $k\in\N$, $f\in L^1(\mu^k)$, and $\eta_{\neq}^k$ is the set of all $k$-tuples of distinct points of $\eta$. In \cite{LachiezeReyPeccati2011,LachiezeReyPeccati2012} and \cite{ReitznerSchulte2011}, Lachi\`eze-Rey and Peccati and Reitzner and Schulte used the bound \eqref{eqn:Peccatietal} for the Wasserstein distance to derive central limit theorems with explicit rates of convergence for such Poisson functionals occurring in stochastic geometry and random graph theory.
Now Theorem \ref{thm:Theorem1} allows us to replace the Wasserstein distance by the Kolmogorov distance without changing the rate of convergence, which means that the inequality \eqref{eqn:dKdW} is not sharp for this class of Poisson functionals.

The main finding of this work, Theorem \ref{thm:Theorem1}, is refined and proven in a different way in the the subsequent paper \cite{EichelsbacherThaele} by Eichelsbacher and Th\"ale. In \cite{LastPeccatiSchulte}, Last, Peccati, and Schulte further simplify the bounds for the normal approximation of Poisson functionals in Kolmogorov distance from Theorem \ref{thm:Theorem1} and \cite{EichelsbacherThaele} and apply them to several problems, where they provide new presumably optimal rates of convergence. Our result for the normal approximation of U-statistics of Poisson point processes (see Theorem \ref{thm:dKUstatistic}) is used by Reitzner, Schulte, and Th\"ale in \cite{ReitznerSchulteThaele2013} and \cite{SchulteThaele2014} to derive central limit theorems with rates of convergence for the Kolmogorov distance for the total edge length of the Gilbert graph and for distances between non-intersecting Poisson $k$-flats, respectively.

This paper is organized in the following way. Before we prove our main result Theorem \ref{thm:Theorem1} in Section \ref{sec:Proof}, we introduce some facts from Malliavin calculus and Stein's method in Section \ref{sec:Preliminaries}. The applications of Theorem \ref{thm:Theorem1} are discussed in Section \ref{sec:Ustatistics}, and the result for U-statistics is shown in Section \ref{sec:ProofUstatistic}.

In this paper, we use the following notation. By $L^p(\P_\eta)$, $p>0$, we denote the set of random variables $Y$ depending on a Poisson point process $\eta$ such that $\E |Y|^p<\infty$. Let $L^p(\mu^n)$, $p>0$, be the set of functions $f:X^n\rightarrow\overline{\R}:=\R\cup\left\{\pm\infty\right\}$ satisfying $\int_{X^n} |f|^p \, \dint \mu^n=\int_{X^n} |f(x_1,\hdots,x_n)|^p \, \mu^n(\dint(x_1,\hdots,x_n))<\infty$ and let $||\cdot||_n$ and $\langle\cdot,\cdot\rangle_{L^2(\mu^n)}$ be the norm and the inner product in $L^2(\mu^n)$, respectively. By $L^p_s(\mu^n)$ we denote the set of all functions $f\in L^p(\mu^n)$ that are symmetric, i.e.\ invariant under permutations of their arguments.

\section{Preliminaries}\label{sec:Preliminaries}

\paragraph{Malliavin calculus for Poisson functionals.} In the sequel, we briefly introduce three Malliavin operators and some properties of them that are necessary for the proofs in this paper. For more details on Malliavin calculus for Poisson functionals we refer to \cite{LP, NualartVives1990, Peccatietal2010,Privault2009} and the references therein.

By $I_n(\cdot)$, $n\geq 1$, we denote the $n$-th multiple Wiener-It\^o integral, which is defined for all functions $f\in L_s^2(\mu^n)$ and satisfies $\E I_n(f)=0$. The multiple Wiener-It\^o integrals are orthogonal in the sense that
\begin{equation}\label{eqn:orthogonality}
\E I_m(f) I_n(g)=\begin{cases} n! \, \langle f,g \rangle_{L^2(\mu^n)}, & m=n\\ 0, & m\neq n\end{cases}
\end{equation}
for all $f\in L_s^2(\mu^m), g\in L^2_s(\mu^n), m,n\geq 1$. We use the convention $I_0(c)=c$ for $c\in\R$. It is known (see \cite{LP} for a proof) that every Poisson functional $F\in L^2(\P_\eta)$ has a unique so-called \textbf{Wiener-It\^o chaos expansion}
\begin{equation}\label{eqn:WienerChaos}
F=\E F+\sum_{n=1}^{\infty}I_n(f_n)
\end{equation}
with $f_n\in L^2_s(\mu^n)$, where the series converges in $L^2(\P_\eta)$. In the following, we call the functions $f_n$ kernels of the Wiener-It\^o chaos expansion of $F$ and say that $F$ has a chaos expansion of order $k$ if $f_n=0$ for all $n>k$. Combining \eqref{eqn:orthogonality} and \eqref{eqn:WienerChaos}, we obtain
$$\V F=\sum_{n=1}^\infty n! \, ||f_n||_n^2.$$
The representation \eqref{eqn:WienerChaos} allows us to define the difference operator $D$, the Ornstein-Uhlenbeck generator $L$, and the Skorohod integral $\delta$ in the following way:

\begin{definition}
Let $F\in L^2(\P_\eta)$ with the Wiener-It\^o chaos expansion \eqref{eqn:WienerChaos}. If\\ $\sum_{n=1}^{\infty}n \, n! \, ||f_n||_n^2<\infty,$
then the random function $z\mapsto D_zF$ defined by
$$D_zF=\sum_{n=1}^{\infty}n \, I_{n-1}(f_n(z,\cdot))$$
is called the \textbf{difference operator} of $F$. For $\sum_{n=1}^{\infty}n^2 \, n! \, ||f_n||_n^2<\infty$ the \textbf{Ornstein-Uhlenbeck generator} of $F$, denoted by $LF$, is given by
$$LF=-\sum_{n=1}^{\infty}n \, I_{n}(f_n).$$
Let $z\mapsto g(z)$ be a random function with a chaos expansion
$$g(z)=g_0(z)+\sum_{n=1}^{\infty}I_{n}(g_n(z,\cdot)), \ g_n(z,\cdot)\in L_s^2(\mu^n),$$
for $\mu$-almost all $z\in X$ and $\sum_{n=0}^{\infty}(n+1)! \, ||g_n||_{n+1}^2<\infty$. Then the \textbf{Skorohod integral} of $g$ is the random variable $\delta(g)$ defined by
$$\delta(g)=\sum_{n=0}^{\infty}I_{n+1}(\tilde{g}_n),$$
where $\tilde{g}_{n}$ is the symmetrization
$\tilde{g}_{n}(x_1,\hdots,x_{n+1})=\frac{1}{(n+1)!}\sum_{\sigma}g_n(x_{\sigma(1)},\hdots,x_{\sigma(n+1)})$ over all permutations $\sigma$ of the $n+1$ variables.
\end{definition}
We denote the domains of these operators by $\dom D$, $\dom  L$, and $\dom \delta$. The difference operator also has the geometric interpretation
\begin{equation}\label{eqn:addone}
D_zF=F(\eta+\delta_z)-F(\eta)
\end{equation}
a.s.\ for $\mu$-almost all $z\in X$, where $\delta_z$ stands for the Dirac measure concentrated at the point $z\in X$, whence it is sometimes called add-one-cost operator (see Theorem 3.3 in \cite{LP}). If $F\notin\dom D$, we can define the difference operator by \eqref{eqn:addone}. If we iterate this definition and put $D^n_{x_1,\hdots,x_n}F=D_{x_n}D^{n-1}_{x_1,\hdots,x_{n-1}}F$, the kernels of the Wiener-It\^o chaos expansion of $F$ in \eqref{eqn:WienerChaos} are given by the formula
$$f_n(x_1,\hdots,x_n)=\frac{1}{n!}\E D^n_{x_1,\hdots,x_n}F =\frac{1}{n!} \E \sum_{I\subset\{1,\hdots,n\}}(-1)^{n+|I|} F(\eta+\sum_{i\in I}\delta_{x_i})$$
for $x_1,\hdots,x_n\in X$ (see Theorem 1.3 in \cite{LP}).

For centred random variables $F\in L^2(\P_\eta)$, i.e.\ $\E F=0$, the inverse Ornstein-Uhlenbeck generator is given by
$$L^{-1}F=-\sum_{n=1}^{\infty}\frac{1}{n}I_n(f_n).$$
The following lemma summarizes how the operators from Malliavin calculus are related.

\begin{lemma}
\begin{itemize}
\item [a)] For every $F\in\dom L$ it holds that $F\in\dom D$, $DF\in\dom \delta$, and
\begin{equation}\label{eqn:IdentityMalliavin}
\delta D F=-LF.
\end{equation}
\item [b)] Let $F\in \dom D$ and $g \in \dom \delta$. Then
\begin{equation}\label{eqn:IntegrationByParts}
\E\langle DF,  g\rangle_{L^2(\mu)}=\E[F \, \delta(g)].
\end{equation}
\end{itemize}
\end{lemma}
For proofs we refer to \cite{Peccatietal2010} and \cite{NualartVives1990}, respectively. Equation \eqref{eqn:IntegrationByParts} is sometimes called  \textbf{integration by parts formula}. Because of this identity, one can see the difference operator and the Skorohod integral as dual operators.

For our applications in Section \ref{sec:Ustatistics} we need a special integration by parts formula, where it is not required that the first Poisson functional is in $\dom D$. In this case, the difference operator is given by \eqref{eqn:addone}.

\begin{lemma}\label{lem:Integration}
Let $F\in L^2(\P_\eta)$, $s\in\R$, and $g\in\dom \delta$ such that $g(z)$ has a Wiener-It\^o chaos expansion of order $k$ for $\mu$-almost all $z\in X$. Moreover, assume that $D_z\1(F>s) \, g(z)\geq 0$ a.s. for $\mu$-almost all $z\in X$. Then
$$\E\langle D\1(F>s),g\rangle_{L^2(\mu)}=\E[\1(F>s) \, \delta(g)].$$
\end{lemma}
\begin{bew}
It is easy to see that $\1(F>s)\in L^2(\P_\eta)$, whence it has a Wiener-It\^o chaos expansion
$$\1(F>s)=\sum_{n=0}^\infty I_n(h_n)$$
with $h_0=\E\1(F>s)$ and kernels $h_n\in L^2_s(\mu^n)$, $n\geq 1$, that are given by
$$h_n(x_1,\hdots,x_n)=\frac{1}{n!} \, \E D^n_{x_1,\hdots,x_n}\1(F>s).$$
For a fixed $z\in X$ the expression $D_z\1(F>s)$ is bounded, and its chaos expansion has the kernels
$$\frac{1}{n!} \, \E D_{x_1,\hdots,x_n}^n D_z\1(F>s)=\frac{1}{n!} \, \E D^{n+1}_{z,x_1,\hdots,x_n}\1(F>s)=(n+1) \, h_{n+1}(z,x_1,\hdots,x_{n}).$$
Hence, we obtain the representation
$$D_z\1(F>s)=\sum_{n=1}^\infty n \, I_{n-1}(h_n(z,\cdot))$$
for all $z\in X$. From Fubini`s theorem and the orthogonality of the multiple Wiener-It\^o integrals it follows that
\begin{equation}\label{eqn:lhsIntegration}
\begin{split}
& \E\langle D\1(F>s),g\rangle_{L^2(\mu)}\\
& = \int_X \E \left[D_z\1(F>s) \, g(z)\right] \, \mu(\dint z)\\
& = \int_{X} \E \left[\sum_{n=1}^\infty n \, I_{n-1}(h_n(z,\cdot)) \sum_{n=0}^k I_n(g_n(z,\cdot))\right] \, \mu(\dint z)\\
& = \int_X \sum_{n=1}^{k+1} n! \int_{X^{n-1}} h_n(z,x_1,\hdots,x_{n-1}) \, g_{n-1}(z,x_1,\hdots,x_{n-1}) \, \mu^{n-1}(\dint (x_1,\hdots,x_{n-1})) \, \mu(\dint z).
\end{split}
\end{equation}
On the other hand, we have
\begin{equation}\label{eqn:rhsIntegration}
 \begin{split}
\E\left[\1(F>s)\,\delta(g)\right]
& = \E\left[\sum_{n=0}^\infty I_n(h_n) \sum_{n=0}^k I_{n+1}(\tilde{g}_n)\right]\\
& = \sum_{n=1}^{k+1} n! \int_{X^n} h_n(x_1,\hdots,x_n) \, \tilde{g}_{n-1}(x_1,\hdots,x_n) \, \mu^n(\dint (x_1,\hdots,x_n))\\
& = \sum_{n=1}^{k+1} n! \int_{X^n} h_n(x_1,\hdots,x_n) \, g_{n-1}(x_1,\hdots,x_n) \, \mu^n(\dint (x_1,\hdots,x_n)).
 \end{split}
\end{equation}
Here, we use the symmetry of $h_n$ in the last step. Comparing \eqref{eqn:lhsIntegration} and \eqref{eqn:rhsIntegration} concludes the proof.
\end{bew}\\
The next lemma provides an upper bound for the second moment of a Skorohod integral, which is used in Section 5.

\begin{lemma}\label{lem:SkorohodSquare}
Let $f\in L^2(\mu^{k+1})$ be symmetric in its last $k$ arguments and let $g(z)=I_k(f(z,\cdot))$. Then
$$\E\left[\delta(g)^2\right]\leq (k+1) \, \E\int_X I_k(f(z,\cdot))^2 \, \mu(\dint z).$$
\end{lemma}
\begin{bew}
By the definition of $\delta$, we obtain $\delta(g)=I_{k+1}(\tilde{f})$ with the symmetrization $$\tilde{f}(x_1,\hdots,x_{k+1})=\frac{1}{(k+1)!}\sum_{\sigma}f(x_{\sigma(1)},\hdots,x_{\sigma(k+1)})$$
as above. From the Cauchy-Schwarz inequality, it follows that $||\tilde{f}||^2_{k+1}\leq ||f||^2_{k+1}$. Combining this with Fubini's Theorem, we have
$$\E\left[\delta(g)^2\right]=(k+1)! \, ||\tilde{f}||^2_{k+1}\leq (k+1)! \, ||f||^2_{k+1}=(k+1) \, \E\int_X I_k(f(z,\cdot))^2 \, \mu(\dint z),$$
which concludes the proof.
\end{bew}

\paragraph{Stein's method.} Besides Malliavin calculus our proof of Theorem \ref{thm:Theorem1} rests upon Stein's method, which goes back to Charles Stein \cite{Stein1972,Stein1986} and is a powerful tool for proving limit theorems. For a detailed and more general introduction into this topic we refer to \cite{ChenGoldsteinShao2011, ChenShao, Stein1986}. Very fundamental for this approach is the following lemma (see Chapter II in \cite{Stein1986}):

\begin{lemma}
For $s\in\R$ the function
\begin{equation}\label{eqn:definitiongt}
g_s(w)=e^{\frac{w^2}{2}}\int_{-\infty}^{w}\left(\1(u\in(-\infty,s])-\P(N\leq s)\right)e^{-\frac{u^2}{2}} \, \dint u
\end{equation}
is a solution of the differential equation
\begin{eqnarray}\label{eqn:Steinsequation}
g_s'(w)-w \, g_s(w)=\1(w\in(-\infty,s])-\P(N\leq s)
\end{eqnarray}
and satisfies
\begin{equation}\label{eqn:gt1}
0<g_s(w)\leq\frac{\sqrt{2\pi}}{4},\ \ |g_s'(w)|\leq 1,\ \ \text{ and }\ \ |w \, g_s(w)|\leq 1
\end{equation}
for any $w\in\R$.
\end{lemma}
Equation \eqref{eqn:Steinsequation} is usually called \textbf{Stein's equation}. The function $g_s$ is infinitely differentiable on $\R\setminus\left\{s\right\}$, but it is not differentiable in $s$. We denote the left-sided and right-sided limits of the derivatives in $s$ by $g_s^{(m)}(s-)$ and $g_s^{(m)}(s+)$, respectively. For the first derivative, a direct computation proves
\begin{equation}\label{eqn:gtt}
g_s'(s+)=-1+g_s'(s-),
\end{equation}
and we define $g_s'(s):=g_s'(s-)$.

By replacing $w$ by a random variable $Z$ and taking the expectation in \eqref{eqn:Steinsequation}, one obtains
$$\E[g_s'(Z)-Z \, g_s(Z)]=\P(Z\leq s)-\P(N\leq s)$$
and as a consequence of the definition of the Kolmogorov distance
\begin{equation}\label{eqn:dKgt}
d_{K}(Z,N)=\sup\limits_{s\in\R}|\E[g_s'(Z)-Z \, g_s(Z)]|.
\end{equation}
The identity \eqref{eqn:dKgt} will be our starting point in Section \ref{sec:Proof}. Note furthermore, that we obtain, by combining \eqref{eqn:Steinsequation} and \eqref{eqn:gt1}, the upper bound
\begin{equation}\label{eqn:gt3}
|g_s''(w)|\leq \frac{\sqrt{2\pi}}{4}+|w|
\end{equation}
for $w\in \R\setminus\left\{s\right\}$.

\section{Proof of Theorem \ref{thm:Theorem1}}\label{sec:Proof}

By a combination of Malliavin calculus and Stein's method similar to that in \cite{Peccatietal2010}, we derive the upper bound for the Kolmogorov distance.
\paragraph{\textbf{Proof of Theorem \ref{thm:Theorem1}}:}Using the identity \eqref{eqn:IdentityMalliavin} and the integration by parts formula \eqref{eqn:IntegrationByParts}, we obtain
\begin{equation}\label{eqn:TrickMalliavin}
\E[F \, g_s(F)]=\E[LL^{-1}F \, g_s(F)]=\E[\delta(-DL^{-1}F) \, g_s(F)]=\E\langle -DL^{-1}F,Dg_s(F)\rangle_{L^2(\mu)}.
\end{equation}
In order to compute $D_zg_s(F)$, we fix $z\in X$ and consider the following cases:
\begin{enumerate}
\item $F,F+D_zF\leq s$ or $F,F+D_zF> s$;
 \item $F\leq s<F+D_zF$;
\item $F+D_zF\leq s<F$.
\end{enumerate}
For $F,F+D_zF\leq s$ or $F,F+D_zF> s$, it follows from \eqref{eqn:addone} and Taylor expansion that
\begin{eqnarray*}
D_z g_s(F) &=&g_s(F+D_zF)-g_s(F)=g_s'(F)D_z F+\frac{1}{2}g_s''(\tilde{F})(D_zF)^2\\ &=:& g_s'(F)D_zF+r_1(F,z,s),
\end{eqnarray*}
where $\tilde{F}$ is between $F$ and $F+D_zF$.
For $F\leq s<F+D_zF$, we obtain by \eqref{eqn:addone}, Taylor expansion, and \eqref{eqn:gtt}
\begin{eqnarray*}
D_z g_s(F) &=& g_s(F+D_zF)-g_s(F)=g_s(F+D_zF)-g_s(s)+g_s(s)-g_s(F)\\ &=& g_s'(s+)(F+D_zF-s)+\frac{1}{2}g_s''(\tilde{F}_1)(F+D_zF-s)^2\\ &&+g_s'(F)(s-F)+\frac{1}{2}g_s''(\tilde{F}_2)(s-F)^2\\
&=& g_s'(F)D_zF+ (g_s'(s-)-1-g_s'(F))(F+D_zF-s)\\ && +\frac{1}{2}g_s''(\tilde{F}_1)(F+D_zF-s)^2+\frac{1}{2}g_s''(\tilde{F}_2)(s-F)^2\\
&=& g_s'(F)D_zF-(F+D_zF-s)+g_s''(\tilde{F}_0)(s-F)(F+D_zF-s)\\ && +\frac{1}{2}g_s''(\tilde{F}_1)(F+D_zF-s)^2+\frac{1}{2}g_s''(\tilde{F}_2)(s-F)^2\\
&=:& g_s'(F)D_zF-(F+D_zF-s)+r_2(F,z,s)
\end{eqnarray*}
with $\tilde{F}_0,\tilde{F}_1,\tilde{F}_2\in \left(F,F+D_zF\right)$. For $F+D_zF\leq s<F$, we have analogously
\begin{eqnarray*}
D_z g_s(F) &=& g_s(F+D_zF)-g_s(F)=g_s(F+D_zF)-g_s(s)+g_s(s)-g_s(F)\\ &=& g_s'(s-)(F+D_zF-s)+\frac{1}{2}g_s''(\tilde{F}_1)(F+D_zF-s)^2\\ &&+g_s'(F)(s-F)+\frac{1}{2}g_s''(\tilde{F}_2)(s-F)^2\\
&=& g_s'(F)D_zF+ (g_s'(s+)+1-g_s'(F))(F+D_zF-s)\\ && +\frac{1}{2}g_s''(\tilde{F}_1)(F+D_zF-s)^2+\frac{1}{2}g_s''(\tilde{F}_2)(s-F)^2\\
&=& g_s'(F)D_zF+(F+D_zF-s)+g_s''(\tilde{F}_0)(s-F)(F+D_zF-s)\\ && +\frac{1}{2}g_s''(\tilde{F}_1)(F+D_zF-s)^2+\frac{1}{2}g_s''(\tilde{F}_2)(s-F)^2\\
&=:& g_s'(F)D_zF+(F+D_zF-s)+r_2(F,z,s)
\end{eqnarray*}
with $\tilde{F}_0,\tilde{F}_1,\tilde{F}_2\in \left(F+D_zF,F\right)$. Thus, $D_zg_s(F)$ has a representation
\begin{equation}\label{eqn:RepresentationDifferenceOperator}
D_zg_s(F)=g_s'(F)D_zF+R(F,z,s),
\end{equation}
where $R(F,z,s)$ is given by
\begin{equation*}
\begin{split}
& R(F,z,s)\\
& =\left(\1(F,F+D_zF\leq s)+\1(F,F+D_zF>s)\right)r_1(F,z,s)\\
& \quad +\left(\1(F\leq s<F+D_zF)+\1(F+D_zF\leq s<F)\right)(r_2(F,z,s)-|F+D_zF-s|).
\end{split}
\end{equation*}
Combining \eqref{eqn:TrickMalliavin} and \eqref{eqn:RepresentationDifferenceOperator} yields
$$\E\left[g_s'(F)-Fg_s(F)\right]=\E\left[g_s'(F)-\langle g_s'(F)DF+R(F,\cdot,s), -DL^{-1}F\rangle_{L^2(\mu)}\right].$$
Thus, the triangle inequality and $|g_s'(F)|\leq 1$ lead to
\begin{eqnarray}\label{eqn:boundStein}
|\E\left[g_s'(F)-Fg_s(F)\right]| & \leq & |\E\left[g_s'(F)\left(1-\langle DF,-DL^{-1}F\rangle_{L^2(\mu)}\right)\right]|\\ &&+|\E\langle R(F,\cdot,s),DL^{-1}F \rangle_{L^2(\mu)}|\notag\\
& \leq & \E |1-\langle DF, -DL^{-1}F\rangle_{L^2(\mu)}|+\E\langle |R(F,\cdot,s)|,|DL^{-1}F|\rangle_{L^2(\mu)}.\notag
\end{eqnarray}
In $r_2(F,z,s)$, we assume that $s$ is between $F$ and $F+D_zF$ so that
$$|F+D_zF-s|\leq |D_zF|\ \text{ and }\ |F-s|\leq |D_zF|.$$
The inequality \eqref{eqn:gt3} and the fact that $\tilde{F}_i$ is between $F$ and $F+D_zF$ allow us to bound all second derivatives in $R(F,z,s)$ by
$$|g''_s(\tilde{F}_i)|\leq\frac{\sqrt{2\pi}}{4}+|F|+|D_zF|.$$
Now it is easy to see that
\begin{equation*}
\begin{split}
& |R(F,z,s)| \\
& \leq \left(\1(F,F+D_zF\leq s)+\1(F,F+D_zF>s)\right)\frac{1}{2}\left(\frac{\sqrt{2\pi}}{4}+|F|+|D_zF|\right)(D_zF)^2\\ & \quad +\left(\1(F\leq s<F+D_zF)+\1(F+D_zF\leq s<F)\right)|D_zF|\\
& \quad +\left(\1(F\leq s<F+D_zF)+\1(F+D_zF\leq s<F)\right)2\left(\frac{\sqrt{2\pi}}{4}+|F|+|D_zF|\right)(D_zF)^2\\
& \leq 2\left(\frac{\sqrt{2\pi}}{4}+|F|+|D_zF|\right)(D_zF)^2\\
& \quad +\left(\1(F\leq s < F+D_zF)+\1(F+D_zF\leq s<F)\right)|D_zF|.
\end{split}
\end{equation*}
By \eqref{eqn:addone}, the last summand can be rewritten as
$$\left(\1(F\leq s < F+D_zF)+\1(F+D_zF\leq s<F)\right)|D_zF|= D_z\1(F>s) \, D_zF.$$
Hence, it follows directly that
\begin{equation*}
\begin{split}
&\E\langle |R(F,\cdot,s)|,|DL^{-1}F|\rangle_{L^2(\mu)}\\ & \leq 2\E\langle (DF)^2,|DL^{-1}F|\rangle_{L^2(\mu)}+2\E\langle (DF)^2,|F \ DL^{-1}F|\rangle_{L^2(\mu)}\\ & \quad +2\E\langle (DF)^2,|DF \ DL^{-1}F|\rangle_{L^2(\mu)}+\E\langle D\1(F>s) \, DF, |DL^{-1}F|\rangle_{L^2(\mu)}.
\end{split}
\end{equation*}
Putting this in \eqref{eqn:boundStein} concludes the proof of the first inequality in \eqref{eqn:BoundNormalApproximation}. The second bound in \eqref{eqn:BoundNormalApproximation} is a direct consequence of the Cauchy-Schwarz inequality.\hfill $\Box$\medskip

In \cite{Peccatietal2010}, the right-hand side of \eqref{eqn:TrickMalliavin} is evaluated for twice differentiable functions $f:\R\rightarrow\R$ with $\sup_{x\in\R}|f'(x)|\leq 1$ and $\sup_{x\in\R}|f''(x)|\leq 2$ (for the Wasserstein distance the solutions of Stein's equation must have these properties) instead of the functions $g_s$ as defined in \eqref{eqn:definitiongt}. For such a function $f$ it holds that
$$D_zf(F)=f'(F) D_zF +\tilde{r}(F)$$
with $|\tilde{r}(F)|\leq (D_zF)^2$. Since this representation is easier than the representation we obtain for $D_zg_s(F)$, the bound for the Wasserstein distance in \eqref{eqn:Peccatietal} is shorter and easier to evaluate than the bound for the Kolmogorov distance in \eqref{eqn:BoundNormalApproximation}.

\section{Applications of Theorem \ref{thm:Theorem1}}\label{sec:Ustatistics}

\paragraph{Normal approximation of a Poisson random variable.} As a first application of Theorem \ref{thm:Theorem1}, we compute an upper bound for the Kolmogorov distance between a Poisson random variable $Y$ with parameter $t>0$ and a normal distribution. In Example 3.5 in \cite{Peccatietal2010}, the bound \eqref{eqn:Peccatietal} is used to compute a bound for the Wasserstein distance, and the known optimal rate of convergence $t^{-\frac{1}{2}}$ is obtained.

$Y$ has the same distribution as $F_t=\sum_{x\in\eta_t}1$, where $\eta_t$ is a Poisson point process on $[0,1]$ with $t$ times the the restriction of the Lebesgue measure on $[0,1]$ as intensity measure $\mu_t$. In the following, we denote by $I_{n,t}(\cdot)$ the $n$-th multiple Wiener-It\^o integral with respect to $\eta_t$. The representation
$$I_{1,t}(f)=\sum_{x\in\eta_t}f(x)-\int_X f(x) \, \mu_t(\dint x)$$
for a Wiener-It\^o integral of a function $f\in L^1(\mu_t)\cap L^2(\mu_t)$ and the fact that
$$F_t=t \int_0^1 1\, \dint x +\sum_{x\in\eta_t} 1 - t \int_0^1 1\, \dint x$$
imply that $F_t$ has the Wiener-It\^o chaos expansion $F_t=\E F_t+I_{1,t}(f_1)=t+I_{1,t}(1)$.
Hence, the standardized random variable
$$G_t=\frac{F_t-\E F_t}{\sqrt{\V F_t}}=\frac{F_t-t}{\sqrt{t}}$$
has the chaos expansion
$G_t=I_{1,t}(1)/\sqrt{t}$
and $D_zG_t=-D_zL^{-1}G_t=1/\sqrt{t}$ for $z\in\left[0,1\right]$. It is easy to see that
$$\E|1-\langle DG_t,-DL^{-1}G_t\rangle_{L^2(\mu_t)}|=|1-\frac{1}{t}\langle1,1\rangle_{L^2(\mu_t)}|=|1-\frac{t}{t}|=0.$$
We obtain
$$\E\langle (DG_t)^2, (DL^{-1}G_t)^2\rangle_{L^2(\mu_t)}=\E\langle (DG_t)^2,(DG_t)^2\rangle_{L^2(\mu_t)}=\frac{1}{t},$$
$\E\langle DG_t,DG_t \rangle_{L^2(\mu_t)}^2=1$, and $\E G_t^4= 3+1/t$ by analogous computations. Since $D_z\1(G_t>s)\,D_zG_t \, |D_zL^{-1}G_t|\geq 0$ for $s\in\R$ and $z\in [0,1]$ and $D_zG_t\,|D_zL^{-1}G_t|=1/t$ for $z\in[0,1]$, it follows from Lemma \ref{lem:Integration} and the Cauchy-Schwarz inequality that
\begin{equation*}
 \begin{split}
 & \sup\limits_{s\in\R}\E\langle D\1(G_t>s), DG_t \, |DL^{-1}G_t|\rangle_{L^2(\mu_t)}\\ & = \sup\limits_{s\in\R} \E[ \1(G_t>s) \, \delta(DG_t \, |DL^{-1}G_t|)] \leq
 \E[\delta(DG_t \, |DL^{-1}G_t|)^2]^{\frac{1}{2}} =\frac{1}{t} \, \E[ I_{1,t}(1)^2]^{\frac{1}{2}} =\frac{1}{\sqrt{t}}.
 \end{split}
\end{equation*}
Now Theorem \ref{thm:Theorem1} yields
$$d_K\left(\frac{Y-t}{\sqrt{t}},N\right)\leq 2 \left(\frac{1}{\sqrt{t}}+\left(3+\frac{1}{t}\right)^{\frac{1}{4}}+1\right)\frac{1}{\sqrt{t}}+\frac{1}{\sqrt{t}}\leq\frac{8}{\sqrt{t}}$$
for $t\geq 1$, which is the classical Berry-Esseen inequality with the optimal rate of convergence (up to a constant).

\paragraph{Normal approximation of U-statistics of Poisson point processes.} As a second application of Theorem \ref{thm:Theorem1}, we discuss U-statistics of the form
$$F=\sum_{(x_1,\hdots,x_k)\in\eta^k_{\neq}}f(x_1,\hdots,x_k)$$
with $k\in\N$ and $f\in L_s^1(\mu^k)$. Here, $\eta^k_{\neq}$ stands for the set of all $k$-tuples of distinct points of $\eta$. If the intensity measure $\mu$ is non-atomic, $\eta$ has no multiple points, and $\eta^k_{\neq}$ can be written as
$$
\eta_{\neq}^k=\left\{(x_1,\hdots,x_k)\in\eta^{k}: x_i\neq x_j \, \forall \, i\neq j\right\}.
$$
In case that $\mu$ has atoms, one has to take into account that distinct points of $\eta$ can have the same location. We denote $k$ as the order of the U-statistic $F$. From now on, we always assume that $F\in L^2(\P_\eta)$ and $\V F>0$. In \cite{ReitznerSchulte2011}, the chaos expansions of such Poisson functionals are investigated, and the bound \eqref{eqn:Peccatietal} is used to prove a central limit theorem with a rate of convergence for the Wasserstein distance. From there it is known that the kernels of the chaos expansion of a U-statistic $F$ are
\begin{equation}\label{eqn:kernelsUstatistic}
f_i(x_1,\hdots,x_i)=\binom{k}{i}\int_{X^{k-i}}f(x_1,\hdots,x_i,y_1,\hdots,y_{k-i}) \, \mu^{k-i}(\dint(y_1,\hdots,y_{k-i}))
\end{equation}
for $i=1,\hdots,k$ and $f_i=0$ for $i>k$. An application of the bound \eqref{eqn:Peccatietal} to such Poisson functionals yields (see Theorem 4.1 in \cite{ReitznerSchulte2011})
\begin{equation}\label{eqn:dWUstatisticBasic}
d_W\left(\frac{F-\E F}{\sqrt{\V F}},N\right)\leq k\sum_{i,j=1}^k\frac{\sqrt{R_{ij}}}{\V F}+k^{\frac{7}{2}}\sum_{i=1}^k\frac{\sqrt{\tilde{R}_i}}{\V F},
\end{equation}
where $R_{ij}$ and $\tilde{R}_i$ are given by
\begin{eqnarray*}
R_{ij}&=&\E \left(\int_X I_{i-1}\left(f_i(z,\cdot)\right)I_{j-1}\left(f_j(z,\cdot)\right) \, \mu(\dint z)\right)^2\\ &&-\left(\E\int_X I_{i-1}\left(f_i(z,\cdot)\right) I_{j-1}\left(f_j(z,\cdot)\right) \, \mu(\dint  z)\right)^2\\
\tilde{R}_i&=&\E\int_X I_{i-1}\left(f_i(z,\cdot)\right)^4 \, \mu(\dint z)
\end{eqnarray*}
for $i,j=1,\hdots,k$. In \cite{ReitznerSchulte2011}, the right-hand side of \eqref{eqn:dWUstatisticBasic} is bounded by a sum of deterministic integrals depending on $f$. Due to technical reasons it is assumed that the U-statistic $F$ is absolutely convergent, which means that the U-statistic $\overline{F}$ given by
$$\overline{F}=\sum_{(x_1,\hdots,x_k)\in\eta^k_{\neq}} |f(x_1,\hdots,x_k)|$$
is in $L^{2}(\P_\eta)$. The U-statistic $\overline{F}$ has a finite Wiener-It\^o chaos expansion with kernels
$$\overline{f}_i(x_1,\hdots,x_i)=\binom{k}{i}\int_{X^{k-i}}|f(x_1,\hdots,x_i,y_1,\hdots,y_{k-i})| \, \mu^{k-i}(\dint(y_1,\hdots,y_{k-i}))$$
for $i=1,\hdots,k$ and $\overline{f}_i=0$ for $i>k$. In order to bound the right-hand side of \eqref{eqn:dWUstatisticBasic} by a sum of deterministic integrals, we use the following notation. For $i,j=1,\hdots,k$ let $\overline{\Pi}_{\geq 2}(i,i,j,j)$ be the set of all partitions $\pi$ of
$$x_1^{(1)},\hdots,x_i^{(1)},x_1^{(2)},\hdots,x_i^{(2)},x_1^{(3)},\hdots,x_j^{(3)},x_1^{(4)},\hdots,x_j^{(4)}$$
such that
\begin{itemize}
\item two variables with the same upper index are in different blocks of $\pi$;
\item each block of $\pi$ includes at least two variables;
\item there are no sets $A_1,A_2\subset \{1,2,3,4\}$ with $A_1\cup A_2=\{1,2,3,4\}$ and $A_1\cap A_2=\emptyset$ such that each block of $\pi$ either consists of variables with upper index in $A_1$ or of variables with upper index in $A_2$.
\end{itemize}
Let $|\pi|$ stand for the number of blocks of a partition $\pi$. For functions $g_1,g_2: X^i\rightarrow\overline{\R}$ and $g_3,g_4: X^j\rightarrow\overline{\R}$ the tensor product $g_1\otimes g_2 \otimes g_3 \otimes g_4: X^{2i+2j}\rightarrow\overline{\R}$ is given by
\begin{equation*}
\begin{split}
& (g_1\otimes g_2 \otimes g_3 \otimes g_4)(x_1^{(1)},\hdots,x_{j}^{(4)})\\ & = g_1(x_1^{(1)},\hdots,x_i^{(1)}) \, g_2(x_1^{(2)},\hdots,x_i^{(2)}) \, g_3(x_1^{(3)},\hdots,x_j^{(3)}) \, g_4(x_1^{(4)},\hdots,x_j^{(4)}).
\end{split}
\end{equation*}
For $\pi\in \overline{\Pi}_{\geq 2}(i,i,j,j)$ we define the function $(g_1\otimes g_2 \otimes g_3 \otimes g_4)_\pi: X^{|\pi|}\rightarrow\overline{\R}$ by replacing all variables that are in the same block of $\pi$ by a new common variable. Since we later integrate over all these new variables, their order does not matter. Using this notation, we define
\begin{eqnarray*}
M_{ij}(f) & = & \sum_{\pi\in \overline{\Pi}_{\geq 2}(i,i,j,j)} \int_{X^{|\pi|}} (\overline{f}_i\otimes\overline{f}_i\otimes \overline{f}_j\otimes \overline{f}_j)_\pi \, \dint\mu^{|\pi|}
\end{eqnarray*}
for $i,j=1,\hdots,k$. Now we can state the following upper bound for the Wasserstein distance (see Theorem 4.7 in \cite{ReitznerSchulte2011}):

\begin{proposition}\label{prop:dWUstatistic}
Let $F\in L^2(\P_\eta)$ be an absolutely convergent U-statistic of order $k$ with $\V F>0$ and let $N$ be a standard Gaussian random variable. Then
\begin{equation}\label{eqn:dWUstatistic}
d_W\left(\frac{F-\E F}{\sqrt{\V F}},N\right)\leq 2k^{\frac{7}{2}}\sum_{1\leq i\leq j\leq k}\frac{\sqrt{M_{ij}(f)}}{\V F}.
\end{equation}
\end{proposition}
In this situation, we can use Theorem \ref{thm:Theorem1} to prove the following bound analogous to \eqref{eqn:dWUstatistic} for the Kolmogorov distance between a standardized U-statistic and a standard Gaussian random variable:

\begin{theorem}\label{thm:dKUstatistic}
Let $F\in L^2(\P_\eta)$ be an absolutely convergent U-statistic of order $k$ with $\V F>0$ and let $N$ be a standard Gaussian random variable. Then
\begin{equation}\label{eqn:dKUstatistic}
d_K\left(\frac{F-\E F}{\sqrt{\V F}},N\right)\leq 19 k^{5} \sum_{i,j=1}^k\frac{\sqrt{M_{ij}(f)}}{\V F}.
\end{equation}
\end{theorem}
Before we prove this theorem in Section \ref{sec:ProofUstatistic}, we discuss some of its consequences. Let us consider a family of Poisson point processes $\eta_t$ with $\sigma$-finite intensity measures $\mu_t$  and U-statistics
$$F_t=\sum_{(x_1,\hdots,x_k)\in(\eta_t)^k_{\neq}}f_t(x_1,\hdots,x_k)$$
with $f_t\in L^1_s(\mu^k_t)$ and $F_t\in L^2(\P_{\eta_t})$ such that
$$\frac{\sqrt{M_{ij}(f_t)}}{\V F_t}\rightarrow 0 \ \text{ as }\ t\rightarrow\infty \ \text{ for all } \ i,j=1,\hdots,k.$$
Here, we integrate with respect to $\mu_t$ in $M_{ij}(f_t)$. Comparing the right-hand sides in \eqref{eqn:dWUstatistic} and \eqref{eqn:dKUstatistic} for the U-statistics $F_t$, we see that the bounds for the Wasserstein and Kolmogorov distance have the same rates of convergence and differ only by constants. An important special case of the described setting is that the Poisson point process depends on a real valued intensity parameter. The following corollary deals with this situation and is the counterpart of Theorem 5.2 in \cite{ReitznerSchulte2011} for the Kolmogorov distance.

\begin{corol}
Let $\eta_t$ be a Poisson point process with intensity measure $\mu_t=t\mu$ with $t\geq 1$ and a fixed $\sigma$-finite measure $\mu$ and let $N$ be a standard Gaussian random variable. We consider U-statistics $F_t\in L^2(\P_{\eta_t})$ of the form
$$F_t=g(t)\sum_{(x_1,\hdots,x_k)\in(\eta_t)^k_{\neq}} f(x_1,\hdots,x_k)$$
with $g: (0,\infty)\rightarrow (0,\infty)$ and $f\in L^1_s(\mu^k)$ independent of $t$. Moreover, we assume that
$$\int_X \left(\int_{X^{k-1}}f(x,y_1,\hdots,y_{k-1}) \, \mu^{k-1}(\dint (y_1,\hdots,y_{k-1}))\right)^2 \, \mu(\dint x)>0$$
and $M_{ij}(f)<\infty$ for $i,j=1,\hdots,k$. Then there is a constant $C>0$ such that
$$d_K\left(\frac{F_t-\E F_t}{\sqrt{\V F_t}},N\right)\leq C \, t^{-\frac{1}{2}}$$
for all $t\geq 1$.
\end{corol}
This corollary follows from bounding $\sqrt{M_{ij}(f_t)}/\V F_t$ in the same way as in the proof of \cite[Theorem 5.2]{ReitznerSchulte2011}.

In \cite{LachiezeReyPeccati2011}, a so-called fourth moment condition for Poisson functionals with positive variance, finite Wiener-It\^o chaos expansion, and non-negative kernels satisfying some integrability conditions is derived. More precisely, for such Poisson functionals it is proven that
$$d_W\left(\frac{F-\E F}{\sqrt{\V F}},N\right)\leq C_{W,k} \sqrt{\E\left(\frac{F-\E F}{\sqrt{\V F}}\right)^4-3}$$
with a constant $C_{W,k}>0$ only depending on $k$.
U-statistics $F\in L^2(\P_\eta)$ of the form
$$F=\sum_{(x_1,\hdots,x_k)\in\eta^k_{\neq}}f(x_1,\hdots,x_k)\ \text{ with }f\in L^1_s(\mu^k)\text{ and }f\geq 0$$
such that $\V F>0$ and $ M_{ij}(f)<\infty$ for $i,j=1,\hdots,k$ belong to this class and satisfy
$$\frac{M_{ij}(f)}{(\V F)^2}\leq \E \left(\frac{F-\E F}{\sqrt{\V F}}\right)^4-3.$$
Then \eqref{eqn:dKUstatistic} can be modified to
$$d_K\left(\frac{F-\E F}{\sqrt{\V F}},N\right)\leq C_k \sqrt{\E\left(\frac{F-\E F}{\sqrt{\V F}}\right)^4-3}$$
with a constant $C_k>0$ only depending on $k$.

\section{Proof of Theorem \ref{thm:dKUstatistic}}\label{sec:ProofUstatistic}

In our proof of Theorem \ref{thm:dKUstatistic}, we make use of the following property of U-statistics:

\begin{lemma}\label{lem:Linverse}
For a U-statistic $F\in L^2(\P_\eta)$ of order $k$ the inverse of the Ornstein-Uhlenbeck generator has a representation
\begin{equation}\label{eqn:Linverse}
\begin{split}
& -L^{-1}(F-\E F)\\
& = \sum_{m=1}^k\frac{1}{m}\sum_{(x_1,\hdots,x_m)\in\eta^m_{\neq}}\int_{X^{k-m}} f(x_1,\hdots,x_m,y_1,\hdots,y_{k-m}) \, \mu^{k-m}(\dint (y_1,\hdots,y_{k-m}))\\ & \quad -\sum_{m=1}^k\frac{1}{m}\int_{X^k} f(y_1,\hdots,y_k) \, \mu^k(\dint (y_1,\hdots,y_k)).
\end{split}
\end{equation}
\end{lemma}

\begin{bew}
We define $\widehat{f}_i:X^i\rightarrow\overline{\R}$ by $\widehat{f}_i(x_1,\hdots,x_i)=\binom{k}{i}^{-1}f_i(x_1,\hdots,x_i)$ for $i=1,\hdots,k$. Using this notation and formula \eqref{eqn:kernelsUstatistic} for the kernels of a U-statistic, we obtain the chaos expansion
\begin{eqnarray*}
&&\sum_{(x_1,\hdots,x_m)\in\eta^m_{\neq}}\int_{X^{k-m}}f(x_1,\hdots,x_m,y_1,\hdots,y_{k-m}) \, \mu^{k-m}(\dint (y_1,\hdots,y_{k-m}))\notag\\ &&=\int_{X^k}f(y_1,\hdots,y_k) \, \mu^k(\dint(y_1,\hdots,y_{k}))+\sum_{i=1}^{m}\binom{m}{i}I_i(\widehat{f}_i)
\end{eqnarray*}
for $m=1,\hdots,k$. Combining this with an identity for binomial coefficients, we see that the right-hand side in \eqref{eqn:Linverse} equals
\begin{eqnarray*}
\sum_{m=1}^k\frac{1}{m}\sum_{i=1}^{m}\binom{m}{i}I_i(\widehat{f}_i)&=&\sum_{m=1}^k\sum_{i=1}^{k}\frac{1}{m}\binom{m}{i}I_i(\widehat{f}_i)=\sum_{i=1}^k\sum_{m=1}^k\frac{1}{m}\binom{m}{i}I_i(\widehat{f}_i)\\ &=&\sum_{i=1}^k\frac{1}{i}\binom{k}{i}I_i(\widehat{f}_i)=\sum_{i=1}^k\frac{1}{i}I_i(f_i),
\end{eqnarray*}
which is the chaos expansion of $-L^{-1}(F-\E F)$ by definition.
\end{bew}

In order to deal with expressions as $R_{ij}$ and $\tilde{R}_i$ in the previous section, one needs to compute the expectation of products of multiple Wiener-It\^o integrals. This can be done by using Proposition 6.1 in \cite{Last2014} (see also \cite[Theorem 3.1]{Sur}, \cite[Proposition 4.5.6]{Privault2009}, \cite[Proposition 6.5.1]{PeccatiTaqqu2011}, or \cite[Theorem 3.1]{LPST2012}). This so-called product formula gives us the Wiener-It\^o chaos expansion of a product of two multiple Wiener-It\^o integrals. Together with \eqref{eqn:orthogonality}, one obtains that the expectation of a product of four multiple Wiener-It\^o integrals is a sum of deterministic integrals depending on the integrands of the stochastic integrals and partitions of their variables as used for the definition of $M_{ij}(f)$.

By using this product formula, one can prove in a similar way as in \cite[Subsection 4.2]{ReitznerSchulte2011} that
\begin{equation}\label{eqn:inequalityRM}
R_{ij}\leq M_{ij}(f)\ \text{ and }\ \int_X \E I_{i-1}(f_i(z,\cdot))^2I_{j-1}(f_j(z,\cdot))^2 \, \mu(\dint z)\leq M_{ij}(f)\
\end{equation}
for $i,j=1,\hdots,k$. Moreover, we prepare the proof of Theorem \ref{thm:dKUstatistic} by the following two lemmas:

\begin{lemma}\label{lem:fourthMoment}
Let $F\in L^2(\P_\eta)$ be an absolutely convergent U-statistic with $M_{ij}(f)<\infty$ for $i,j=1,\hdots,k$. Then
$$\E (F-\E F)^4 \leq  k^2\sum_{i,j=1}^k M_{ij}(f)+3k^2 (\V F)^2.$$
\end{lemma}
\begin{bew}
Using the Wiener-It\^o chaos expansion of $F$ and the Cauchy-Schwarz inequality, we obtain
$$\E (F-\E F)^4 = \E \left(\sum_{i=1}^k I_i(f_i)\right)^2 \left(\sum_{j=1}^k I_j(f_j)\right)^2 \leq k^2 \, \E \sum_{i=1}^k I_i(f_i)^2 \sum_{j=1}^k I_j(f_j)^2.$$
Now the previously mentioned product formula for multiple Wiener-It\^o integrals and $\V F=\sum_{n=1}^k n! \, \|f_n\|_n^2$ yield that
$$\E I_i(f_i)^2I_j(f_j)^2 \leq M_{ij}(f)+ 3 \, i! \, ||f_i||_i^2 \, j! \, ||f_j||_j^2 \leq M_{ij}(f)+ 3 (\V F)^2.$$
In the first inequality, we have equality for $i=j$ and $f\geq 0$.
\end{bew}

\begin{lemma}\label{lem:sup}
Let $F\in L^2(\P_\eta)$ be an absolutely convergent U-statistic with $M_{ij}(f)<\infty$ for $i,j=1,\hdots,k$. Then
\begin{equation*}
\begin{split}
& \sup\limits_{s\in\R}\E\langle D\1(F>s),DF \, |DL^{-1}(F-\E F)|\rangle_{L^2(\mu)}\\
& \leq \sqrt{(2k-1) \, \E\langle (D\overline{F})^2, (DL^{-1}\left(\overline{F}-\E\overline{F}\right))^2\rangle_{L^2(\mu)}}.
\end{split}
\end{equation*}
\end{lemma}
\begin{bew}
We can write the U-statistic $F$ as
$$F=\sum_{(x_1,\hdots,x_k)\in\eta^k_{\neq}}f(x_1,\hdots,x_k)=\underbrace{\sum_{(x_1,\hdots,x_k)\in\eta^k_{\neq}}f^{+}(x_1,\hdots,x_k)}_{=F^{+}}-\underbrace{\sum_{(x_1,\hdots,x_k)\in\eta^k_{\neq}}f^{-}(x_1,\hdots,x_k)}_{=F^{-}}$$
with $f^{+}=\max\left\{f,0\right\}$ and $f^{-}=\max\left\{-f,0\right\}$ and have $\overline{F}=F^++F^-$. As a consequence of \eqref{eqn:addone}, we know that $D_zV\geq 0$ for a U-statistic $V$ with non-negative summands. Combining this with $f^+,f^-\geq 0$ and Lemma \ref{lem:Linverse}, we see that
$$-D_zL^{-1}\left(F^{+}-\E F^{+}\right)\geq 0\quad \text{ and }\quad -D_zL^{-1}\left(F^{-}-\E F^{-}\right)\geq 0.$$
Moreover, it holds that $D_z\1(F>s) \, D_zF\geq 0$. Proposition 6.1 in \cite{Last2014} implies that the product $D_zF\ D_zL^{-1}\left(\overline{F}-\E\overline{F}\right)$ has a finite chaos expansion with an order less than or equal to $2k-2$. Together with Lemma \ref{lem:Integration}, the Cauchy-Schwarz inequality, and Lemma \ref{lem:SkorohodSquare}, we obtain
\begin{equation*}
\begin{split}
&\sup\limits_{s\in\R}\E\langle D\1(F>s),DF \, |DL^{-1}(F-\E F)|\rangle_{L^2(\mu)}\\
&= \sup\limits_{s\in\R}\E\langle D\1(F>s),DF \, |DL^{-1}\left(F^{+}-\E F^{+}-F^{-}+\E F^{-}\right)|\rangle_{L^2(\mu)}\\
&\leq \sup\limits_{s\in\R}\E\langle D\1(F>s),DF \left(-DL^{-1}(F^{+}-\E F^{+})-DL^{-1}(F^{-}-\E F^{-})\right)\rangle_{L^2(\mu)}\\
&\leq \E\left[\delta\left(DF\, DL^{-1}\left(\overline{F}-\E\overline{F}\right)\right)^2\right]^{\frac{1}{2}}\\
& \leq \sqrt{(2k-1)\E\langle (DF)^2, (DL^{-1}\left(\overline{F}-\E\overline{F}\right))^2\rangle_{L^2(\mu)}}.
\end{split}
\end{equation*}
Now the fact that $(D_zF)^2\leq (D_z\overline{F})^2$ concludes the proof.
\end{bew}

\paragraph{\textbf{Proof of Theorem \ref{thm:dKUstatistic}}:}
In the following, we can assume that $M_{ij}(f)<\infty$ for $i,j=1,\hdots,k$ since \eqref{eqn:dKUstatistic} is obviously true, otherwise. We consider the standardized random variable
$G=(F-\E F)/\sqrt{\V F}$, whose Wiener-It\^o chaos expansion has the kernels $g_i=f_i/\sqrt{\V F}$ for $i=1,\hdots,k$. In order to simplify our notation, we use the abbreviation
$$S=\sum_{i,j=1}^k \frac{\sqrt{M_{ij}(f)}}{\V F}.$$
Exactly as in the proof of Theorem 4.1 in \cite{ReitznerSchulte2011}, we obtain
$$\E\left|1-\langle DG,-DL^{-1}G\rangle_{L^2(\mu)}\right|\leq k\sum_{i,j=1}^k\frac{\sqrt{R_{ij}}}{\V F}\leq k\sum_{i,j=1}^k\frac{\sqrt{M_{ij}(f)}}{\V F} = k S.$$
From straightforward computations using Fubini's Theorem, the Cauchy-Schwarz inequality, and \eqref{eqn:inequalityRM}, it follows that
\begin{eqnarray*}
\E\langle (DG)^2,(DL^{-1}G)^2\rangle_{L^2(\mu)}&=&\int_X \E (D_zG)^2(D_zL^{-1}G)^2 \, \mu(\dint  z)\\
&=& \int_X \E \left(\sum_{i=1}^{k} i \, I_{i-1}(g_i(z,\cdot))\right)^2\left(\sum_{i=1}^{k} I_{i-1}(g_i(z,\cdot))\right)^2 \, \mu(\dint z)\\
&\leq & k^4 \int_X \E \sum_{i=1}^{k} I_{i-1}(g_i(z,\cdot))^2 \sum_{j=1}^{k} I_{j-1}(g_j(z,\cdot))^2 \, \mu(\dint z)\\
& \leq & k^4 \sum_{i,j=1}^k \frac{M_{ij}(f)}{(\V F)^2} \leq k^4 S^2,
\end{eqnarray*}
\begin{eqnarray*}
\E\langle (DG)^2,(DG)^2\rangle_{L^2(\mu)} &=& \int_X \E\left(\sum_{i=1}^{k}i \, I_{i-1}(g_i(z,\cdot))\right)^4 \, \mu(\dint z)\\
& \leq & k^6 \int_X \E \sum_{i=1}^k I_{i-1}(g_i(z,\cdot))^2 \sum_{j=1}^k I_{j-1}(g_j(z,\cdot))^2 \, \mu(\dint z)\\
& \leq & k^6 \sum_{i,j=1}^k \frac{M_{ij}(f)}{(\V F)^2}\leq k^6 S^2,
\end{eqnarray*}
and
\begin{eqnarray*}
\E\langle DG,DG\rangle_{L^2(\mu)}^2 &=& \E\left(\sum_{i,j=1}^k i j\int_X I_{i-1}(g_i(z,\cdot)) I_{j-1}(g_j(z,\cdot)) \, \mu(\dint z)\right)^2\\
&\leq & k^2 \sum_{i,j=1}^k i^2 j^2 \E\left(\int_X I_{i-1}(g_i(z,\cdot)) I_{j-1}(g_j(z,\cdot)) \, \mu(\dint z)\right)^2\\ &\leq & k^6 \sum_{i,j=1}^k \frac{R_{ij}}{(\V F)^2}+k^4\sum_{i=1}^k \frac{(i!)^2||f_i||_i^4}{(\V F)^2}\\
&\leq& k^6 \sum_{i,j=1}^k \frac{R_{ij}}{(\V F)^2}+k^4\leq k^6 \sum_{i,j=1}^k \frac{M_{ij}(f)}{(\V F)^2}+k^4\leq k^6S^2+k^4.
\end{eqnarray*}
As a consequence of Lemma \ref{lem:fourthMoment}, we have that
$$\E\left(\frac{F-\E F}{\sqrt{\V F}}\right)^4\leq k^2\sum_{i,j=1}^k\frac{M_{ij}(f)}{(\V F)^2}+3 k^2\leq k^2 S^2 + 3k^2.$$
Combining the last three inequalities, we obtain
\begin{equation*}
 \begin{split}
  & 2\sqrt{\E\langle (DG)^2,(DG)^2\rangle_{L^2(\mu)}}+2\left(\E\langle DG, DG\rangle_{L^2(\mu)}^2\right)^{\frac{1}{4}} \left((E G^4)^{\frac{1}{4}}+1\right)\\
  & \leq 2 k^3 S+2(k^{\frac{3}{2}}\sqrt{S}+k)(\sqrt{k}\sqrt{S}+3^{\frac{1}{4}}\sqrt{k}+1)\leq 16k^3
 \end{split}
\end{equation*}
for $S\leq 1$.

Lemma \ref{lem:sup} together with a similar computation as for $\E\langle (D G)^2, (DL^{-1}G)^2\rangle_{L^2(\mu)}$ implies that
$$\sup_{s\in\R}\E\langle D\1(G>s),DG \, |DL^{-1}G|\rangle_{L^2(\mu)} \leq \sqrt{(2k-1)k^4\sum_{i,j=1}^k\frac{M_{ij}(f)}{(\V F)^2}}\leq \sqrt{2}k^{\frac{5}{2}} S.$$
Thus, it follows from Theorem \ref{thm:Theorem1} that
$$d_K\left(\frac{F-\E F}{\sqrt{\V F}},N\right) \leq k S+16 k^3 \, k^2 S +\sqrt{2} k^{\frac{5}{2}} S\leq 19k^5 S$$
for $S\leq 1$. Otherwise, this inequality still holds since the Kolmogorov distance is at most one.
\hfill $\Box$\medskip

In a similar way, one can also obtain an upper bound for the Kolmogorov distance between a Gaussian random variable and a finite sum of Poisson U-statistics. This class of Poisson functionals is interesting since Theorem 3.6 in \cite{ReitznerSchulte2011} tells us that each Poisson functional $F\in L^2(\P_\eta)$ of order $k$ with kernels $f_i\in L^1_s(\mu^i)\cap L_s^2(\mu^i)$ for $i=1,\hdots,k$ is a finite sum of Poisson U-statistics (and a constant). For such a Poisson functional the upper bounds for the inner products in the proof of Theorem \ref{thm:dKUstatistic} that depend on $R_{ij}$ and
$$\int_X \E I_{i-1}(f_i(z,\cdot))^2I_{j-1}(f_j(z,\cdot))^2\mu(\dint z)$$
for $i,j=1,\hdots,k$ still hold. Moreover, we can compute a similar bound as in Lemma \ref{lem:sup} using the representation of $F$ as a sum of Poisson U-statistics. Together with the fourth centred moment of $F$, one can obtain an upper bound for the Kolmogorov distance between $(F-\E F)/\sqrt{\V F}$ and a standard Gaussian random variable.

\vspace{1cm}
{\bf Acknowledgement:} Large parts of this paper were written during a stay at Case Western Reserve University (February to July 2012) supported by the German Academic Exchange Service. The author thanks Elizabeth Meckes, Giovanni Peccati, Matthias Reitzner, and Christoph Th\"ale for some valuable hints and helpful discussions.


\begin{thebibliography}{30}\small

\bibitem{ChenGoldsteinShao2011}
Chen, L., Goldstein, L., Shao, Q.: Normal Approximation by Stein's Method. Springer, Berlin (2011)

\bibitem{ChenShao}
Chen, L., Shao, Q.: Stein's method for normal approximation. An Introduction to Stein's Method, 1--59, Singapore University Press, Singapore (2005)

\bibitem{EichelsbacherThaele}
Eichelsbacher, P., Th\"ale, C.: New Berry-Esseen bounds for non-linear functionals of Poisson random measures. arXiv: 1310.1595 (2013)

\bibitem{LachiezeReyPeccati2011}
Lachi\`eze-Rey, R., Peccati, G.: Fine Gaussian fluctuations on the Poisson space, I:
contractions, cumulants and geometric random graphs. Electron. J. Probab. \textbf{18}, Article 32 (2013)

\bibitem{LachiezeReyPeccati2012}
Lachi\`eze-Rey, R., Peccati, G.: Fine Gaussian fluctuations on the Poisson space, II:
rescaled kernels, marked processes and geometric U-statistics. Stochastic Process. Appl. \textbf{123}, 4186--4218 (2013)

\bibitem{Last2014}
Last, G.: Stochastic analysis for Poisson processes. arXiv: 1405.4416 (2014)

\bibitem{LastPeccatiSchulte}
Last, G., Peccati, G., Schulte, M.: Normal approximation on Poisson spaces: Mehler's formula, second order Poincar\'e inequalitites and stabilization. arXiv: 1401.7568 (2014)

\bibitem{LP}
Last, G., Penrose, M.D.: Poisson process Fock space representation, chaos expansion and covariance inequalities. Probab. Theory Related Fields \textbf{150}, 663--690 (2011)

\bibitem{LPST2012}
Last, G., Penrose, M.D., Schulte, M., Th\"ale, C.: Moments and central limit theorems for some multivariate Poisson functionals. Adv. in Appl. Probab. \textbf{46}, 348--364 (2014)

\bibitem{NualartVives1990}
Nualart, D., Vives, J.: Anticipative calculus for the Poisson process based on the Fock space. Lecture Notes in Math. \textbf{1426}, 154--165 (1990)

\bibitem{Peccatietal2010}
Peccati, G., Sol\'e, J.L., Taqqu, M.S., Utzet, F.: Stein's method and normal approximation of Poisson functionals. Ann. Probab. \textbf{38}, 443--478 (2010)

\bibitem{PeccatiTaqqu2011}
Peccati, G., Taqqu, M.S.: Wiener Chaos: Moments, Cumulants and Diagram Formulae. Springer, Berlin (2011)

\bibitem{Privault2009}
Privault, N.: Stochastic Analysis in Discrete and Continuous Settings with Normal Martingales. Springer, Berlin (2009)

\bibitem{ReitznerSchulte2011}
Reitzner, M., Schulte, M.: Central limit theorems for U-statistics of Poisson point processes. Ann. Probab. \textbf{41}, 3879--3909 (2013)

\bibitem{ReitznerSchulteThaele2013}
Reitzner, M., Schulte, M., Th\"ale, C.: Limit theory for the Gilbert graph. arXiv: 1312.4861 (2013)

\bibitem{SchulteThaele2014}
Schulte, M., Th\"ale, C.: Distances between Poisson $k$-flats. Methodol. Comput. Appl. Probab. \textbf{16}, 311--329 (2014)

\bibitem{Stein1972}
Stein, C.: A bound for the error in the normal approximation to the distribution of a sum of dependent random variables. Proceedings of the Sixth Berkley Symposium on Mathematical Statistics Probability 2, Univ. California Press, Berkley (1972)

\bibitem{Stein1986}
Stein, C.: Approximate Computation of Expectations. Institute of Mathematical Statistics, Hayward, CA (1986)

\bibitem{Sur}
Surgailis, D.: On multiple Poisson stochastic integrals and associated Markov semigroups. Probab. Math. Statist. \textbf{38}, 217--239 (1984)

\end{thebibliography}
\end{document}